\documentclass{article} 
\usepackage{iclr2024_conference, times}
 \usepackage{amsmath,amsfonts}

\usepackage{multicol}
 \usepackage{graphicx}
 \usepackage{multirow,array,booktabs}
\begin{document}
\title{Physics Informed Modeling of Ecosystem Respiration  via Dynamic Mode Decomposition with Control Input } 
\author{Maha Shadaydeh \\Department of Mathematics and Computer Science, Friedrich Schiller University Jena,  Germany\\
Email: maha.shadaydeh@uni-jena.de
\AND  Joachim Denzler \\Department of Mathematics and Computer Science, Friedrich Schiller University Jena,  Germany
\AND Mirco Migliavacca \\
European Commission Joint Research Centre, 1050 Brussels, Belgium}


\maketitle
\begin{abstract}
Ecosystem respiration (Reco) represents a major component of the global carbon cycle, and accurate characterization of its dynamics is essential for a comprehensive understanding of ecosystem-climate interactions and the impacts of climate extremes on the ecosystem. This paper presents a novel data-driven and physics-aware method for estimating Reco dynamics using the dynamic mode decomposition with control input (DMDc) technique, an emerging tool for analyzing nonlinear dynamical systems. The proposed model represents Reco as a state space model with an autonomous component and an exogenous control input. The control input can be any ecosystem driver(s) such air temperature, soil temperature, or soil water content. This unique modeling approach allows controlled intervention to study the effects of different inputs on the system. Experimental results using Fluxnet2015  data show that the prediction accuracy of Reco dynamics achieved with DMDc is comparable to state-of-the-art methods, making it a promising tool for analyzing the dynamic behavior of different vegetation ecosystems on multi-temporal scales in response to different climatic drivers.  
\end{abstract}
%
\section{Introduction}
\label{intro}
Carbon losses from ecosystems affect climate change \citep{Melillo2017}. Ecosystem respiration (Reco), the sum of autotrophs (plants) and heterotrophs (bacteria, fungi, and animals) respiration, represents a major component of the global carbon cycle. Accurate estimation of Reco dynamics is necessary for a better understanding of ecosystem-climate interactions and climate extremes' impact on ecosystems. This paper presents a data-driven yet physics-aware method for estimating Reco dynamics using the dynamic mode decomposition with control input technique \citep{Proctor2016}, an emerging tool for analyzing nonlinear dynamical systems.

Ecosystem respiration is typically described as an exponential function of temperature based on the law of thermodynamics \citep{Lloyd1994}. However, this exponential relationship is defined on the basis of some parameters, such as temperature sensitivity, which are assumed to remain constant. Several studies have pointed to the dependence of these parameters on other drivers of ecosystems \citep{mahecha}, which could lead to either over- or underestimation of Reco depending on the temperature range. Regression models partially compensate for such problems by using temporally moving windows to estimate the parameters of the model \citep{Reichstein2012}.

The eddy covariance (EC) technique is widely used to measure the net ecosystem exchange (NEE), which is the difference between the total CO2 release due to all respiration processes and the gross carbon uptake by photosynthesis (GPP) \citep{FLUXNET2001}. The two CO2 fluxes, Reco and GPP, are derived using what is known as NEE partitioning methods  \citet{Lloyd1994, Reichestein2005,fluxnet2015}. These approaches use NEE measurements to fit physiologically based nonlinear relationships to estimate GPP and Reco using meteorological drivers. These functional relationships are either a light-response-based function linking global incoming radiation and GPP \citet{Gilmanov} or respiration functions based on temperature \citep{Reichstein2012}. These methods are adopted as standard processing tools in the FLUXNET community \citet{fluxnet2015} and are referred to as daytime (DT) and nighttime (NT) methods, respectively.

Recently, the capabilities of deep neural networks in learning complex and nonlinear dynamical models have been used for modeling Reco dynamics \citep{Tramontana2020, Trifunov2021}. These approaches provide data-driven equation-free estimates of Reco with the flexibility to include other meteorological and biological drivers affecting Reco.  
Despite their improved performance,  these methods require sufficient training data and experimental tuning of the hyper-parameters of the deep networks used. Accordingly, the trained model can not consider short-term variations in ecosystem respiration. 

Koopman operator \citep{Koopman1931} enables the transformation of finite-dimensional nonlinear system dynamics to an infinite-dimensional linear dynamical system. 
However, finding the Koopman operator's eigenfunctions remains a significant obstacle to its implementation. 
Dynamic mode decomposition (DMD) is a simple numerical algorithm approximating the Koopman operator with a best-fit linear model that advances measurements from one time step to the next \citep{Mezi2005, Schmid2010, Rowley2009, Brunton2017}. It is an equation-free system identification method where the underlying dynamics of the system are learned from snap-shots in time of measurement data. DMD decomposes system dynamics into temporal modes whereby each mode represents a fixed oscillation frequency and decay/growth rate. 
These modes are not orthonormal modes and thus less parsimonious than those generated by PCA but more physically meaningful. The extension of DMD to dynamical systems with control inputs is the dynamic mode decomposition with control (DMDc) \citep{Proctor2016}. DMDc provides a robust framework for analyzing the behavior of complex dynamical systems under the influence of external inputs and has been successfully applied in diverse nonlinear dynamical systems, e.g.,  fluids dynamics and neuroscience. 

In this paper, we extend our initial work  \citep{Shadaydeh} and present a novel data-driven  Reco estimation approach that can also serve as an NEE partitioning method using DMDc in a sliding temporal window.  The system state Reco is presented as a state space dynamical model with an autonomous component and an exogenous component that serves as a control input function. The control input to the system could be air temperature (Tair) per the thermodynamics-based exponential function and/or other observed drivers such as soil temperature, soil water contents (SWC), etc.  Such modeling of Reco allows disentangling the exogenous effect of the control input, e.g., Tair,  from the autonomous dynamics of Reco. Subsequently, it facilitates intervention studies with different control inputs to study their effect on the ecosystem. To deal with unobserved drivers of Reco other than temperature or any used control input, we use time delay embedding (TDE) of the history of the system's state and control input. According to the Takens Theory \citep{Takens1981}, embedding the state history guarantees under certain conditions that the system will learn the trajectories of the original system. The model training is performed on night-time NEE of Fluxnet2015 data \citep{fluxnet2015}.  Reco day- and nighttime values, within forecast intervals that span one to several days, are then predicted using the trained model with the future values of the control input(s). The performance of the DMDc with different control inputs, namely Tair and SWC, is compared to DT and NT NEE partitioning methods reported in the Fluxnet data. The impact of TDE on the performance of the proposed DMDc approach is also investigated for long-term forecasting. 

The paper is structured as follows: The Methodology section provides a brief overview of  DMDc and how it is utilized in the analysis of Reco dynamics. It also discusses using TDE and the singular value decomposition of the resulting Hankel matrix of nighttime NEE to characterize the Reco system dynamics. Experimental results using DMDc with different control inputs and delay embedding dimensions are compared to state-of-the-art methods in  Section 3. Section 4 concludes the paper.


\section{Methodology}

\subsection{Dynamic Mode Decomposition with Control (DMDc)}

DMD \citep{Schmid2010} analyzes the relationship between pairs of measurements from a dynamical system. The measurements, $\mathbf{x}(k)$ and $\mathbf{x}(k+1)$, where $k$ indicates the temporal iteration from a discrete dynamical system, are assumed to be approximately related by a linear operator $\mathbf{A} $ such that the relation   $\mathbf{x}(k+1) \approx \mathbf{A} \mathbf{x}(k)$  holds for all pairs of measurements $k$,  where $\mathbf{x} \in \mathbb{R}^{n}$ and $\mathbf{A} \in \mathbb{R}^{n \times n}$. Let $\mathbf{X}=\left[\begin{array}{cccc} \mathbf{x}(1) & \mathbf{x}(2) & \cdots & \mathbf{x}(M-1)\end{array}\right]$ and the time-shifted matrix $\mathbf{X}^{\prime}=\left[\begin{array}{cccc}
\mathbf{x}(2) & \mathbf{x}(3) & \ldots & \mathbf{x}(M) \\
\end{array}\right]$ where $M$ is the total number of snapshots. This relationship can be described in a  matrix form as
   $\mathbf{X}^{\prime} \approx \mathbf{A} \mathbf{X}$.
The best-fit linear operator is computed as  $\mathbf{A}=\mathbf{X}^{\prime} \mathbf{X}^{\dagger}$  where ${ }^{\dagger}$ is the Moore-Penrose pseudoinverse which is computed using the singular value decomposition (SVD) of  $\mathbf{X}$. The eigen-decomposition of $\mathbf{A}$ defined by $\mathbf{A} \mathbf{W}=\mathbf{W} \mathbf{\Lambda}$ yields eigenvalues and eigenvectors that can be used for future state prediction. 
Specifically,  let $\phi_j$ and $\lambda_j$ be the jth eigenvector-eigenvalue pair of $\mathbf{A}$,  $\mathbf{x}(k)$ at any time $k$  is defined as 
\begin{equation}
      \mathbf{x}(k)= \sum_{j=1}^{N} b_j \phi_j \lambda_j^k  
\end{equation}

where $b_j$ denotes the initial state $\mathbf{x}(1)$. 


DMDc \citep{Proctor2016} extends the DMD to dynamical systems with control inputs. The goal of DMDc is to analyze the relationship between a future system measurement $\mathbf{x}(k+1)$ with the current measurement $\mathbf{x}(k)$ and the current control input $\mathbf{u}(k)$ such that $\mathbf{x}(k+1) \approx \mathbf{A} \mathbf{x}(k)+\mathbf{B} \mathbf{u}(k)$,   where $\mathbf{x}(k) \in \mathbb{R}^{n}, \mathbf{u}(k) \in \mathbb{R}^{l}, \mathbf{A} \in \mathbb{R}^{n \times n}$, and $\mathbf{B} \in \mathbb{R}^{n \times l}$. Let  $\mathbf{\Upsilon}=\left[\begin{array}{cccc}
 \mathbf{u}(1) & \mathbf{u}(2) & \cdots & \mathbf{u}(m-1) \\ 
 \end{array}\right]$.   Utilizing the three data matrices, DMDc is focused on finding best-fit approximations to the mappings $\mathbf{A}$ and $\mathbf{B}$ such that 
 \begin{equation}
 \label{Eq: DMDc}
   \mathbf{X}^{\prime} \approx \mathbf{A X}+\mathbf{B \Upsilon} 
 \end{equation}
holds for all trios of data. The augmented snapshot matrix $[\mathbf{X}; \mathbf{\Upsilon}]$ is decomposed using SVD to extract the dynamic modes of the system, including the effects of the control inputs. The DMD algorithm explained above is then applied to the augmented snapshot matrix to identify the dynamic modes of the system. Once the dynamic modes are identified, they can be used to reconstruct the evolution of the system under the influence of control inputs.  The detailed steps of DMDc are summarised in Table \ref{tab:DMDc} in the Appendix.


 \subsection{Modeling the dynamics of Reco using DMDc}
 
 In Eq.  \ref{Eq: DMDc},  let the state of the system $\mathbf{X}$  and the control input $\mathbf{\Upsilon}$  be two vectors which constitute of $M$  nights measurement of nighttime NEE (NEEnight) and the control input, e.g., Tair,  respectively. Specifically, for Tair as a control input, we define 
 \begin{eqnarray}
 \mathbf{X}&=& \left[\begin{array}{cccc}  \rm{NEEnight}(t_1) & \rm{NEEnight}(t_2) & \cdots & \rm{NEEnight}(t_{M-1})\end{array}\right] \nonumber\\ 
 \mathbf{\Upsilon} &=&\left[\begin{array}{cccc} \rm{Tair}(t_1) & \rm{Tair}(t_2) & \cdots & \rm{Tair}(t_{M-1})\end{array}\right].  
 \end{eqnarray}
Here, NEEnight(t$_i), i=1, \cdots, M$ refers to a vector of nighttime measurement of the day $t_i$. The matrices $\mathbf{A}$ and $\mathbf{B}$ in Eq. \ref{Eq: DMDc} are identified using the steps of the DMDc algorithm in Table \ref{tab:DMDc}. 

The dynamical model of Reco is learned using the DMDc method applied to the values of NEEnight and the selected control input, e.g., Tair, using a temporal sliding window of size $M$. The trained model is then used to forecast day and night time values of Reco using the day and night time values of the control input within the forecast horizon, which might span one to several days. 

\subsection{Time delay embedding (TDE) for characterization of Reco dynamics}
 Time delay embedding is a classical approach for dealing with partial state information. Most often, we can observe only partially the drivers affecting the state of the system under observation, i.e., the ecosystem respiration system in this study.    This method augments a state $\mathbf{x}$  with its history.  According to  Takens theory \citep{Takens1981}, under mild conditions on the observables, the dynamics of the augmented state are guaranteed to be diffeomorphic to the dynamics of the original state $\mathbf{x}$, provided that the embedding dimension N satisfies  $N \geq 2n + 1$, where $n$ is the dimension of the state. One method to obtain the TDE dimension N is by computing the SVD of a Hankel matrix \citep{Brunton2017}. The Hankel matrix is formed by stacking the columns of the snapshot matrices $\mathbf{X}$ with delay embedding, e.g., the Hankel matrix of NEEnight in Eq.  \ref{eq:hankel}.  
 
 \begin{eqnarray}
  {\bf{H}}{\rm{ = }}\left[ {\begin{array}{*{20}{c}} \rm{NEEnight}(t_1) & \rm{NEEnight}(t_2) & \cdots & \rm{NEEnight}(t_{M-N}) \\ \\ \rm{NEEnight}(t_2) & \rm{NEEnight}(t_3) & \cdots & \rm{NEEnight}(t_{M-N+1)} \\ \\ \vdots & \vdots & \ddots & \vdots \\ \\ \rm{NEEnight}(t_N) & \rm{NEEnight}(t_{N+1}) & \cdots & \rm{NEEnight}(t_{M}) \\ \end{array}} \right]
  \label{eq:hankel}
 \end{eqnarray}

 The singular values of the  Hankel matrix capture the energy or variance associated with each of the DMDc modes \citep{Brunton2017}. Thus, the number of dominant modes necessary for Reco forecast can be estimated using the rank of the Hankel matrix of NEEnight. 
 
 Analyzing the spread of singular values of the Hankel matrix also provides valuable insights into the underlying system's dynamics.  A well-structured spread of singular values indicates a well-observable or predictable system.  Smaller singular values indicate significant noise in the system.  Comparing the spread of singular values of different sites can help characterize the differences in Reco dynamics in different vegetation and average temperature sites, as shown in the following section. 

 Another advantage of using TDE in the proposed method is that it allows for learning Reco dynamics from short data. Specifically, the delay embedding compensates for using advanced measurement in time. This is relevant as it enables ecosystem forecasts to consider very transient factors affecting short-term variations in the system dynamics, such as perspiration.

\section{Experimental Results and Discussion}
 We used the half-hourly EC Fluxnet2015 data \citep{fluxnet2015}  measured at Fluxnet sites with different vegetation types and average temperature to investigate the performance of our method and the impact of different control inputs, e.g.,  Tair and SWC,  on Reco dynamics. In particular, we selected five different forests of different types: Deciduous Broadleaf Forests (DBF), Evergreen Needleleaf Forests (ENF), Mixed Forests (MF), as well as three Grassland  (GRA) sites. 
 
Since the direct measurement of daytime Reco is not available, the root mean square error (RMSE) between the Reco forecast values and NEEnight during the nights of the forecast interval is used as a validation measure. We used the variable 'Night' of the Fluxnet data as a reference for nighttime. The method is applied on a temporal sliding window over the last two years of several sites of  Fluxnet 2015  data from May to September (January to April for the southern hemisphere). Only good quality data with at least $80\%$ of quality flag $Q<2$ is used. Minimax normalization is applied to control input variables only. We compare the performance of the proposed method using DMDc with and without control input (DMD) to NT and DT Reco values of the  Fluxnet data.   

\begin{figure}[t!]
   \centering
    \includegraphics[width=1\linewidth]{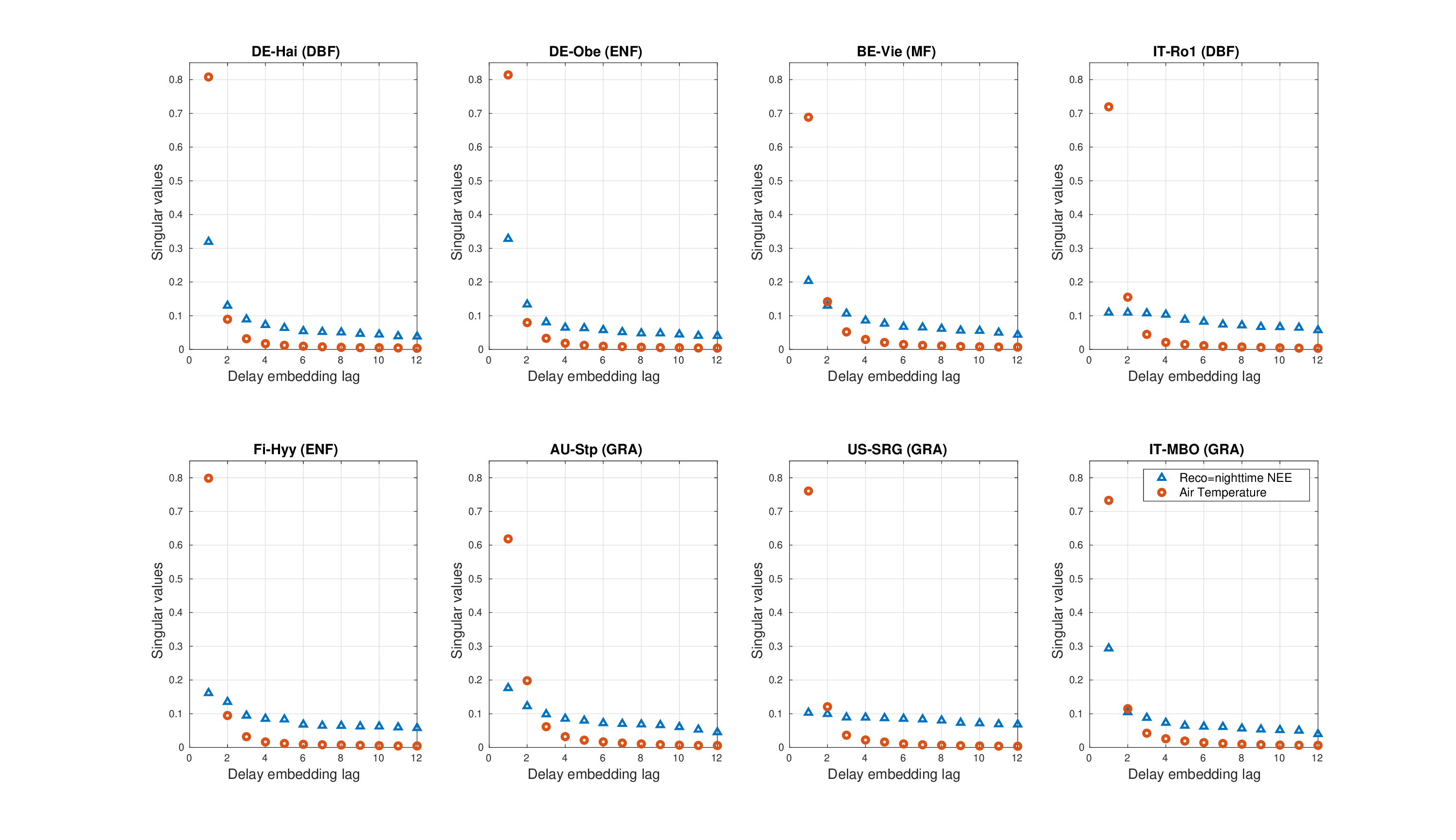}
   \caption{Singular values of the Hankel matrices H of NEEnight and the control input Tair in different Fluxnet sites and vegetation types. The spread of the singular values is an indicator of the predictability of the system. The number of distinct singular values is the number of dominant dynamic modes in the ecosystem. We can observe four to six dominant modes depending on the vegetation type and location. }
    \label{fig:svd}
\end{figure}

Figure \ref{fig:svd} shows the singular values of the Hankel matrices of  NEEnight and Tair for different vegetation sites. In most sites, we can observe four to six dominant modes. The number of modes used in the reconstruction and forecast of Reco is selected for each site accordingly. Results of RMSE of one day forecast using five nights for training are given in Table \ref{tab1:RMSE1}. These results indicate a comparable performance of the proposed method to DT and NT methods. The spread of the singular values of the NEEnight Hankel matrix of the different sites shown in  Figure \ref{fig:svd}  can explain the prediction accuracy of our method.  We can notice that sites with well-spread singular values, such as DE-Hai (DBF) and BE-Vie(MF), are more predictable, and hence our proposed DMDc-based method outperformed DT and NT methods. However, for sites with less structured singular values, such as IT-Ro1(GRA) and AU-Stp (GRA), we can observe a slightly lower accuracy of DMDc. This is expected since the low spread of singular values makes the system theoretically less observable or predictable.  The low spread of singular values could result from significant noise in the system or slow climate dynamics in the site location.  It should be noted, however, that the DMD method on its own, without any control input, is still able to produce reasonable accuracy, which is comparable to the DT and NT methods in some sites.

Our results also show a promising performance in capturing the influence of different ecosystem drivers. We can notice that SWC is a better predictor of Reco dynamics for grasslands or higher average temperature sites, while Tair is better for colder average-temperature forest sites.

\begin{table}[t!]
  \caption{ \small Comparison of different methods’ nighttime Reco forecasts to the ground truth NEEnight. We show the NEEnight estimation RMSE error metric for all methods and vegetation types.   The performance of DT and NT NEE partitioning methods are compared to those of  DMD and DMDc when using five nights of NEEnight for training and one-day Reco forecast, i.e., validation on the next night. The best results are shown in bold; lower is better. } 
    \centering
    \scalebox{0.85}{
    \begin{tabular}{l    l c c c c c c c}             
   \hline 
      & \multicolumn{8}{c} {Fluxnet Site Name (Vegetation Type)} \\ 
\hline 
  Method  & DE-Hai & DE-Obe & BE-Vie  & IT-Ro1 & FY-Hyy &  AU-Stp & US-SRG & IT-MBO \\
  (Control Input)     & (DBF)& (ENF)& (MF)& (DBF) &(ENF) &   (GRA) &  (GRA)& (GRA)\\        
\hline
Fluxnet  NT Reco    & 0.192  &  \textbf{0.190}  &  0.143   & 0.147  &  0.204 &    0.096 &   0.163 &   0.127\\    
\hline
Fluxnet DT Reco      &   0.251  &  0.209  &  0.173   & \textbf{0.126}   & 0.221   & \textbf{0.091}   & \textbf{0.122}  &  0.165\\    
\hline
 DMD    &   0.195  &  0.222  &  0.140  &  0.144  &  0.191   & 0.098  &  0.174 &   0.130\\    
\hline 
DMDc(Tair)   &  \textbf{0.182}  &  0.207  &  \textbf{0.137}  &  0.145  &  \textbf{0.188}   & 0.099  &  0.175  &  \textbf{0.126}\\ \hline 
DMDc(SWC)    & 0.191  &  0.213 &   0.141  &  0.145   & 0.197 &   0.098  &  0.171   & 0.127\\      
     \hline
    \end{tabular}}
    \label{tab1:RMSE1}
\end{table}

\begin{table}[t!]
  \caption{ \small Comparison of different methods’ nighttime Reco forecasts to the ground truth NEEnight. We show the NEEnight estimation RMSE  metric for all methods and vegetation types.  The performance of  DMD and DMDc with different control inputs and embedding dimensions are compared using two weeks of training data and two weeks of Reco forecast, i.e., validation on the nights of the next two weeks.  The best results are shown in bold; lower is better.} 
    \centering
    \scalebox{0.9}{
    \begin{tabular}{l    l c c c c c c c}             
   \hline 
      & \multicolumn{8}{c} {Fluxnet Site Name (Vegetation Type)} \\ 
\hline 
  Method (Control Input,  & DE-Hai & DE-Obe & BE-Vie  & IT-Ro1 & FY-Hyy &  AU-Stp & US-SRG & IT-MBO \\
   Embedding dimension)     & (DBF)& (ENF)& (MF)& (DBF) &(ENF) &   (GRA) &  (GRA)& (GRA)\\        
\hline
 DMD    &    0.216  &  0.3231   & 0.199  &  0.193  &  0.279  &  \textbf{0.114}  &  0.242  &  0.171\\    
\hline 
    DMDc(Tair)   &  0.235 &   0.301 &   \textbf{0.189}  &  \textbf{0.186}  &  0.274  &  0.115  &  0.243 &   \textbf{0.161}\\
  DMDc-TDE(Tair, N=4)  &   0.230 &   0.297 &   0.191  &  0.188 &   0.271 &   0.115  &  0.240   & 0.162\\
  DMDc-TDE(Tair, N=6)   & 0.231 &   \textbf{0.294}   & 0.190   & 0.190  &  \textbf{0.260}  &  0.117  &  0.236  &  0.163\\  
   \hline 
   DMDc(SWC)    &   0.219 &   0.326  &  0.198  &  0.195  &  0.296 &   0.115  &  0.230&    0.180\\
  DMDc-TDE(SWC, N=4)  & 0.216   & 0.316  &  0.194  &  0.204  &  0.281  &  0.116  &  0.225  &  0.180\\
  DMDc-TDE(SWC, N=6)  & 0.219  &  0.315  &  0.198  &  0.192 &   0.279   & 0.122 &  \textbf{0.219} &   0.181\\
  \hline

    \end{tabular}}
    \label{tab1:RMSE2}
\end{table}
In a second set of experiments, we aim to investigate the impact of using time delay embedding to improve the accuracy of longer Reco forecast intervals. To this end, we train the DMDc state-space model with two weeks of NEEnight measurement data and forecast the next two weeks.   Table \ref{tab1:RMSE2} shows the results of DMD and DMDc for two different control inputs, Tair and SWC, and two TDE dimensions, N=4 and 6. These two embedding values are selected based on the number of distinct singular values we can observe in Figure   \ref{fig:svd}. Results indicate that, on the one hand, Tair is a better long-term predictor in most studied sites. On the other hand, We can also observe that TDE enhanced the prediction accuracy with SWC as control input for all studied sites. Surprisingly, DMD with no control input has the highest two-week prediction accuracy in Hai-DE (DBE)  and  AU-Stp (GRA). 
\section{Conclusions}
We proposed a data-driven model for analyzing ecosystem respiration dynamics using DMDc. Such modeling of Reco allows disentangling the exogenous effect of the control input, e.g., Tair, from the autonomous dynamics of Reco. Subsequently, it facilitates interventional studies on the control input to study its effect on the ecosystem. The performance of the DMDc with different control inputs, namely Tair and SWC, is compared to DT and NT NEE partitioning methods reported in the Fluxnet data for different vegetation sites. Our results show a promising performance of our approach in the characterization of Reco dynamics and in capturing the influence of different ecosystem drivers, making it a promising tool for analyzing the dynamic behavior of different vegetation ecosystems in response to climate change. Experimental results also showed that for most sites, time delay embedding of different ecosystem drivers could improve long-term prediction accuracy. Current research efforts are focused on implementing the proposed approach on multitemporal scales for improving its long-term forecast accuracy.

\section{Acknowledgment} 
This research is funded by the German Research Foundation (DFG) grant SH 1682/1-1 and  the Carl Zeiss Foundation within the scope of the program line “Breakthroughs: Exploring Intelligent Systems” for “Digitization — explore the basics, use applications.”


\begin{thebibliography}{18}
\providecommand{\natexlab}[1]{#1}
\providecommand{\url}[1]{\texttt{#1}}
\expandafter\ifx\csname urlstyle\endcsname\relax
  \providecommand{\doi}[1]{doi: #1}\else
  \providecommand{\doi}{doi: \begingroup \urlstyle{rm}\Url}\fi

\bibitem[Baldocchi et~al.(2001)Baldocchi, Falge, Gu, et~al.]{FLUXNET2001}
D.~Baldocchi, E.~Falge, L.~H. Gu, et~al.
\newblock Fluxnet: A new tool to study the temporal and spatial variability of ecosystem-scale carbon dioxide, water vapor, and energy flux densities.
\newblock \emph{Bulletin of the American Meteorological Society}, 82\penalty0 (11):\penalty0 2415--2434, 2001.

\bibitem[{Brunton} et~al.(2017){Brunton}, {Brunton}, {Proctor}, {Kaiser}, and {Kutz}]{Brunton2017}
Steven~L. {Brunton}, Bingni~W. {Brunton}, Joshua~L. {Proctor}, Eurika {Kaiser}, and J.~Nathan {Kutz}.
\newblock {Chaos as an intermittently forced linear system}.
\newblock \emph{Nature Communications}, 8:\penalty0 19, May 2017.
\newblock \doi{10.1038/s41467-017-00030-8}.

\bibitem[Gilmanov et~al.(2003)Gilmanov, Verma, Sims, Meyers, Bradford, Burba, and Suyker]{Gilmanov}
Tagir~G. Gilmanov, Shashi~B. Verma, Phillip~L. Sims, Tilden~P. Meyers, James~A. Bradford, George~G. Burba, and Andrew~E. Suyker.
\newblock Gross primary production and light response parameters of four southern plains ecosystems estimated using long-term co2-flux tower measurements.
\newblock \emph{Global Biogeochemical Cycles}, 17\penalty0 (2), 2003.
\newblock \doi{https://doi.org/10.1029/2002GB002023}.

\bibitem[Koopman(1931)]{Koopman1931}
B.~O. Koopman.
\newblock Hamiltonian systems and transformation in hilbert space.
\newblock 17\penalty0 (5):\penalty0 315--318, 1931.
\newblock ISSN 0027-8424.
\newblock \doi{10.1073/pnas.17.5.315}.

\bibitem[Lloyd \& Taylor(1994)Lloyd and Taylor]{Lloyd1994}
J.~Lloyd and J.~A. Taylor.
\newblock On the temperature dependence of soil respiration.
\newblock \emph{Functional Ecology}, 8\penalty0 (3):\penalty0 315--323, 1994.
\newblock ISSN 02698463, 13652435.

\bibitem[Mahecha et~al.(2010)Mahecha, Reichstein, Carvalhais, Lasslop, Lange, Seneviratne, Vargas, Ammann, Arain, Cescatti, et~al.]{mahecha}
Miguel~D Mahecha, Markus Reichstein, Nuno Carvalhais, Gitta Lasslop, Holger Lange, Sonia~I Seneviratne, Rodrigo Vargas, Christof Ammann, M~Altaf Arain, Alessandro Cescatti, et~al.
\newblock Global convergence in the temperature sensitivity of respiration at ecosystem level.
\newblock \emph{Science}, 329\penalty0 (5993):\penalty0 838--840, 2010.

\bibitem[Melillo et~al.(2017)Melillo, Frey, DeAngelis, Werner, Bernard, Bowles, Pold, Knorr, and Grandy]{Melillo2017}
J.~M. Melillo, S.~D. Frey, K.~M. DeAngelis, W.~J. Werner, M.~J. Bernard, F.~P. Bowles, G.~Pold, M.~A. Knorr, and A.~S. Grandy.
\newblock Long-term pattern and magnitude of soil carbon feedback to the climate system in a warming world.
\newblock \emph{Science}, 358\penalty0 (6359):\penalty0 101--105, 2017.
\newblock \doi{10.1126/science.aan2874}.

\bibitem[Mezi{\'c}(2005)]{Mezi2005}
Igor Mezi{\'c}.
\newblock Spectral properties of dynamical systems, model reduction and decompositions.
\newblock \emph{Nonlinear Dynamics}, 41:\penalty0 309--325, 2005.

\bibitem[Pastorello(2020)]{fluxnet2015}
Gilberto et~al Pastorello.
\newblock The {FLUXNET2015} dataset and the {ONEFlux} processing pipeline for eddy covariance data.
\newblock \emph{Scientific Data}, 7\penalty0 (1):\penalty0 225, 2020.
\newblock ISSN 2052-4463.
\newblock \doi{10.1038/s41597-020-0534-3}.

\bibitem[Proctor et~al.(2016)Proctor, Brunton, and Kutz]{Proctor2016}
Joshua~L. Proctor, Steven~L. Brunton, and J.~Nathan Kutz.
\newblock Dynamic mode decomposition with control.
\newblock \emph{SIAM Journal on Applied Dynamical Systems}, 15\penalty0 (1):\penalty0 142--161, 2016.
\newblock \doi{10.1137/15M1013857}.

\bibitem[Reichstein et~al.(2005)Reichstein, Falge, Baldocchi, et~al.]{Reichestein2005}
Markus Reichstein, Eva Falge, Dennis Baldocchi, et~al.
\newblock On the separation of net ecosystem exchange into assimilation and ecosystem respiration: review and improved algorithm.
\newblock \emph{Global Change Biology}, 11\penalty0 (9):\penalty0 1424--1439, 2005.

\bibitem[Reichstein et~al.(2012)Reichstein, Stoy, Desal, et~al.]{Reichstein2012}
Markus Reichstein, Paul~C. Stoy, Ankur~R. Desal, et~al.
\newblock Partitioning of net fluxes.
\newblock In Marc Aubinet, Timo Vesala, and Dario Papale (eds.), \emph{Eddy Covariance: A Practical Guide to Measurement and Data Analysis}, pp.\  263--289. Springer Netherlands, 2012.
\newblock ISBN 978-94-007-2350-4.

\bibitem[Rowley et~al.(2009)Rowley, Mezic, Bagheri, Schlatter, and Henningson]{Rowley2009}
Clarence~W. Rowley, Igor Mezic, Shervin Bagheri, Philipp Schlatter, and Dan~S. Henningson.
\newblock Spectral analysis of nonlinear flows.
\newblock \emph{Journal of Fluid Mechanics}, 641:\penalty0 115–127, 2009.
\newblock \doi{10.1017/S0022112009992059}.

\bibitem[Schmid \& Peters(2010)Schmid and Peters]{Schmid2010}
Schmid and Peters.
\newblock Dynamic mode decomposition of numerical and experimental data.
\newblock \emph{Journal of Fluid Mechanics}, 656:\penalty0 5–28, 2010.
\newblock \doi{10.1017/S0022112010001217}.

\bibitem[Shadaydeh et~al.(2022)Shadaydeh, Denzler, and Migliavacca]{Shadaydeh}
Maha Shadaydeh, Joachim Denzler, and Mirco Migliavacca.
\newblock Partitioning of net ecosystem exchange using dynamic mode decomposition and time delay embedding.
\newblock \emph{Engineering Proceedings}, 18\penalty0 (1), 2022.
\newblock \doi{10.3390/engproc2022018013}.

\bibitem[Takens(1981)]{Takens1981}
Floris Takens.
\newblock Detecting strange attractors in turbulence.
\newblock In David Rand and Lai-Sang Young (eds.), \emph{Dynamical Systems and Turbulence, Warwick 1980}, pp.\  366--381, Berlin, Heidelberg, 1981. Springer Berlin Heidelberg.
\newblock ISBN 978-3-540-38945-3.

\bibitem[Tramontana et~al.(2020)Tramontana, Migliavacca, Jung, et~al.]{Tramontana2020}
Gianluca Tramontana, Mirco Migliavacca, Martin Jung, et~al.
\newblock Partitioning net carbon dioxide fluxes into photosynthesis and respiration using neural networks.
\newblock \emph{Global Change Biology}, 26\penalty0 (9):\penalty0 52355--5253, 2020.

\bibitem[Trifunov et~al.(2021)Trifunov, Shadaydeh, Runge, Reichstein, and Denzler]{Trifunov2021}
Violeta~Teodora Trifunov, Maha Shadaydeh, Jakob Runge, Markus Reichstein, and Joachim Denzler.
\newblock A data-driven approach to partitioning net ecosystem exchange using a deep state space model.
\newblock \emph{IEEE Access}, 9:\penalty0 107873--107883, 2021.
\newblock \doi{10.1109/ACCESS.2021.3101129}.

\end{thebibliography}
\newpage
\appendix
\section{Appendix}
\begin{table}[th]
\centering
\caption{Main steps of the DMDc Method 
\cite{Proctor2016} }
 \scalebox{1}{
\begin{tabular}{l}
\hline
1.     Collect the system measurement and control snapshots \\
2. Form the matrices $\mathbf{X}, \mathbf{X}^{\prime}$, and $\boldsymbol{\Upsilon}$\\
  3.  Stack the data matrices $\mathbf{X}$ and $\Upsilon$ to construct the matrix $\boldsymbol{\Omega}$ \\
4. Compute the $\mathrm{SVD}$ of $\boldsymbol{\Omega}$ thereby obtaining the decomposition \\
$\boldsymbol{\Omega} \approx \tilde{\mathbf{U}} \tilde{\boldsymbol{\Sigma}} \tilde{\mathbf{V}}^{*}$ with truncation value $p$. \\
5. Compute the SVD of $\mathbf{X}^{\prime}$, thereby obtaining the decomposition \\
$\mathbf{X}^{\prime} \approx \hat{\mathbf{U}} \hat{\mathbf{\Sigma}} \hat{\mathbf{V}}^{*}$ with truncation value $r$ \\
6. Compute the following:\\
\vbox{\begin{eqnarray}
    \tilde{\mathbf{A}}=\hat{\mathbf{U}}^{*} \mathbf{X}^{\prime} \tilde{\mathbf{V}} \tilde{\Sigma}^{-1} \tilde{\mathbf{U}}_{1}^{*} \hat{\mathbf{U}}  \\
   \tilde{\mathbf{B}}=\hat{\mathbf{U}}^{*} \mathbf{X}^{\prime} \tilde{\mathbf{V}} \tilde{\Sigma}^{-1} \tilde{\mathbf{U}}_{2}^{*}\\
  \tilde{\mathbf{A}}=\hat{\mathbf{U}}^{*} \mathbf{X}^{\prime} \tilde{\mathbf{V}} \tilde{\Sigma}^{-1} \tilde{\mathbf{U}}_{1}^{*} \hat{\mathbf{U}} \\
   \tilde{\mathbf{B}}=\hat{\mathbf{U}}^{*} \mathbf{X}^{\prime} \tilde{\mathbf{V}} \tilde{\Sigma}^{-1} \tilde{\mathbf{U}}_{2}^{*}
  \end{eqnarray}}\\
7. Perform the eigenvalue decomposition \\
 \vbox{ \begin{eqnarray}
  \tilde{\mathbf{A}} \mathbf{W}=\mathbf{W} \mathbf{\Lambda} 
 \end{eqnarray}}\\
8. Compute the dynamic modes of the operator $\mathbf{A}$.\\
 \vbox{ \begin{eqnarray}
  \boldsymbol{\Phi}=\mathbf{X}^{\prime} \tilde{\mathbf{V}} \tilde{\boldsymbol{\Sigma}}^{-1} \tilde{\mathbf{U}}_{1}^{*} \hat{\mathbf{U}} \mathbf{W}
\end{eqnarray}}\\
\hline
\end{tabular}}
     \label{tab:DMDc}
 \end{table}

\end{document}